\documentclass[12pt,reqno,b5paper]{amsproc}%
\usepackage{color,anysize,float,CJK}
\usepackage{citeref}
\usepackage[noadjust]{cite}
\usepackage{txfonts,esint,bm}
\usepackage[colorlinks,urlcolor=blue,citecolor=blue,linkcolor=blue]{hyperref}

\marginsize{15mm}{15mm}{1cm}{1cm}
\renewcommand{\bibitempages}[1]{}

\newtheorem{thm}{Theorem}[section]

\newtheorem{lem}[thm]{Lemma}

\theoremstyle{definition}

\theoremstyle{remark}

\newtheorem*{pf}{\it Proof}

\numberwithin{equation}{section}
\allowdisplaybreaks

\begin{document}

\title{Infinitely many solutions for Kirchhoff equations with indefinite potential}
\author[S. Jiang \and S. Liu]{Shuai Jiang\,$^\mathrm{a}$, Shibo Liu\,$^\mathrm{b}$\vspace{-1em}}
\dedicatory{$^\mathrm{a}$School of Mathematical Sciences, Xiamen University\\
Xiamen 361006, P.R. China\\
$^\mathrm{b}$Department of Mathematical Sciences, Florida Institute of Technology\\
Melbourne, FL 32901, USA}
\maketitle

\begin{abstract}
We obtain a sequence of solutions converging to zero for the Kirchhoff
equation%
\[
-\left(  1+\int_{\Omega}\left\vert \nabla u\right\vert ^{2}\right)  \Delta
u+V(x)u=f(u)\text{,\qquad}u\in H_{0}^{1}(\Omega)
\]
via truncating technique and a variant of Clark's theorem due to Liu--Wang,
where $\Omega$ is a bounded smooth domain $\Omega\subset\mathbb{R}^{N}$.
Similar result for Schr\"{o}dinger-Poisson system on a bounded smooth domain
$\Omega\subset\mathbb{R}^{3}$ is also presented.

\end{abstract}

\section{Introduction}

In a recent paper \cite{MR4162412}, He and Wu studied the following elliptic
boundary value problem%
\[
-\Delta u+V(x)u=f(x,u)\text{,\qquad}u\in H_{0}^{1}(\Omega)
\]
with indefinite linear part $-\Delta+V$, where $\Omega\subset\mathbb{R}^{N}$
is a bounded smooth domain and the odd nonlinearity $f:\Omega\times
\mathbb{R}\to\mathbb{R}$ is sublinear at zero:%
\[
\lim_{\left\vert t\right\vert \rightarrow0}\frac{1}{t^{2}}\int_{0}%
^{t}f(x,s)\,\mathrm{d}s=+\infty\text{.}%
\]
Using truncating technique and Liu--Wang's variant of Clark's theorem
\cite[Theorem 1.1]{MR3400440}, they obtained a sequence of solutions
conversing to zero in $H_{0}^{1}(\Omega)$.

Motivated by \cite{MR4162412}, in this note we consider the following
Kirchhoff equation on a bounded smooth domain $\Omega\subset\mathbb{R}^{N}$,%
\begin{equation}
-\left(  1+\int_{\Omega}\left\vert \nabla u\right\vert ^{2}\right)  \Delta
u+V(x)u=f(x,u)\text{,\qquad}u\in H_{0}^{1}(\Omega)\text{.}\label{k1}%
\end{equation}
We impose the following conditions on the potential $V$ and the nonlinearity
$f$,

\begin{enumerate}
\item[$\left(  V\right)  $] $V\in C(\Omega)$ is bounded;

\item[$\left(  f_{1}\right)  $] $f\in C(\Omega\times\mathbb{R})$ is
subcritical, that is%
\[
\lim_{\left\vert t\right\vert \rightarrow\infty}\frac{f(x,t)t}{\left\vert
t\right\vert ^{2^{\ast}}}=0\text{,\quad where }2^{\ast}=\frac{2N}{N-2}\text{
is the critical exponent;}%
\]

\item[$\left(  f_{2}\right)  $] $f(x,\cdot)$ is odd for all $x\in\Omega$,
$f(x,0)=0$, and is sublinear at zero:
\begin{equation}
\lim_{\left\vert t\right\vert \rightarrow0}\frac{F(x,t)}{t^{2}}=+\infty
\text{,\qquad where }F(x,t)=\int_{0}^{t}f(x,s)\,\mathrm{d}s\text{.}\label{k0}%
\end{equation}

\end{enumerate}
We will prove the following theorem.

\begin{thm}
\label{kt1}Suppose $\left(  V\right)  $, $\left(  f_{1}\right)  $ and $\left(
f_{2}\right)  $ hold, then the problem (\ref{k1}) possesses a sequence of
nontrivial solutions converging to zero.
\end{thm}

Boundary value problems of the form (\ref{k1}) are closely related to the wave
equation%
\[
\psi_{tt}-\left(  a+b\int_{\Omega}\left\vert \nabla\psi\right\vert
^{2}\right)  \Delta\psi=g(x,\psi)\text{,\qquad}\left(  t,x\right)  \in\left(
0,T\right)  \times\Omega\text{,}%
\]
which was used by G. Kirchhoff to investigate vibrations of elastic strings
with changing length. Starting from Alves \emph{et al.}\ \cite{MR2123187}, where a
variational approach is developed to solve (\ref{k1}), many existence results
for (\ref{k1}) appear. For example, Cheng \emph{et al.}\  \cite{MR2984595} considered
the case that $V(x)=0$ and the nonlinearity is of the form%
\begin{equation}
f(x,t)=\alpha(x)\left\vert t\right\vert ^{q-2}t+g(x,t)\text{,}\label{f1}%
\end{equation}
where $q\in\left(  1,2\right)  $, $g(x,t)=o(\left\vert t\right\vert )$ as
$t\rightarrow0$. Obviously such $f$ satisfies our assumption $\left(
f_{2}\right)  $. Since they need $H_{0}^{1}(\Omega)\hookrightarrow
L^{r}(\Omega)$ for $r>4$, it is assumed in \cite{MR2984595} that $N\leq3$.
Furtado and Zanata \cite{MR3668404} also considered (\ref{k1}) with $V(x)=0$
and $f$ as in (\ref{f1}); but they only imposed local conditions to $g(x,t)$
for $\left\vert t\right\vert $ small ($g$ needs not be odd and subcritical for
$\left\vert t\right\vert $ large). Using some idea from Wang \cite{MR1828946},
they got a sequence of solutions $\left\{  u_{k}\right\}  $ for the truncated
problem with an odd and subcritical $\tilde{g}$ in place of $g$, $\tilde
{g}(x,t)=g(x,t)$ for $\left\vert t\right\vert $ small; then applied
$L^{\infty}$-estimate to show that $\left\vert u_{k}\right\vert _{\infty
}\rightarrow0$ and concluded that for $k$ large $u_{k}$ are solutions of the original
problem. Since our problem (\ref{k1}) may be indefinite, such $L^{\infty}%
$-estimate seems not applicable, this is why we need $f$ to be globally odd
and subcritical. For more recent papers on Kirchhoff equations, the reader is
referred to \cite{MR4098134,MR4201645,MR4361869}.

When $N=3$, for the following Schr\"{o}dinger-Poisson system on a bounded
smooth domain $\Omega$%
\begin{equation}
\left\{
\begin{array}
[c]{ll}%
-\Delta u+V(x)u+\phi u=f(x,u) & \text{in }\Omega\text{,}\\
-\Delta\phi=u^{2} & \text{in }\Omega\text{,}\\
u=\phi=0 & \text{on }\partial\Omega\text{,}%
\end{array}
\right.  \label{kk}%
\end{equation}
we have similar result.

\begin{thm}
\label{2k}Suppose $\left(  V\right)  $, $\left(  f_{1}\right)  $ and $\left(
f_{2}\right)  $ hold, then the problem (\ref{kk}) possesses a sequence of
nontrivial solutions $\left(  u_{n},\phi_{n}\right)  \rightarrow\left(
0,0\right)  $ in $H_{0}^{1}(\Omega)\times H_{0}^{1}(\Omega)$.
\end{thm}

Since the seminar work or Benci \emph{et al.}\  \cite{MR1659454},
Schr\"{o}dinger-Poisson system has been an active field of research, for
recent work on Schr\"{o}dinger-Poisson system on bounded domain we mention
\cite{MR3279523,MR4208981,MR4299007}.

\section{Proof of Theorem \ref{kt1}}

The dependence on $x$ in $f(x,u)$ is not essential in our discussion of
(\ref{k1}) and (\ref{kk}). Therefore in what follows we write $f(u)$ for
$f(x,u)$ for simplicity.

It is well known that to find weak solutions of (\ref{k1}), it suffices to
find critical points of the $C^{1}$-functional $\Phi:H_{0}^{1}(\Omega
)\rightarrow\mathbb{R}$ defined by%
\begin{equation}
\Phi(u)=\frac{1}{2}\int\left(  \left\vert \nabla u\right\vert ^{2}%
+V(x)u\right)  +\frac{1}{4}\left(  \int\left\vert \nabla u\right\vert
^{2}\right)  ^{2}-\int F(u)\text{,}\label{k2}%
\end{equation}
here and below the integrals are taken over $\Omega$. Let $E^{-}$, $E^{0}$,
and $E^{+}$ be the negative space, null space, and positive space of the
quadratic form (the first term) in (\ref{k2}). For $u\in E:=H_{0}^{1}(\Omega)$, we always
denote by $u^{\pm}$ and $u^{0}$ the orthogonal projections of $u$ on $E^{\pm}$
and $E^{0}$. Because of the condition $\left(  V\right)  $, there is an
equivalent norm $\Vert\cdot\Vert$ on $E$ such that%
\[
\Phi(u)=\frac{1}{2}\left(  \Vert u^{+}\Vert^{2}-\left\Vert u^{-}\right\Vert
^{2}\right)  +\frac{1}{4}\left(  \int\left\vert \nabla u\right\vert
^{2}\right)  ^{2}-\int F(u)\text{.}%
\]
We denote by $\left(  \cdot,\cdot\right)  $ the corresponding inner product.

To prove Theorem \ref{kt1} it suffices to find a sequence of critical points
of $\Phi$. For this purpose, we need the following variant of the Clark's
theorem due to Liu--Wang \cite{MR3400440}.

\begin{thm}
[{\cite[Theorem 1.1]{MR3400440}}]\label{kt2}Let $E$ be a Banach space and
$\Phi\in C^{1}(E,\mathbb{R})$ be an even coercive functional satisfying the
$\left(  PS\right)  $ condition and $\Phi(0)=0$. If for any $k\in\mathbb{N}$,
there is an $k$-dimensional subspace $X_{k}$ and $\rho_{k}>0$ such that%
\begin{equation}
\sup_{X_{k}\cap S_{\rho_{_{k}}}}\Phi<0\text{,} \label{k4}%
\end{equation}
where $S_{r}=\left\{  u\in E|\,\left\Vert u\right\Vert =r\right\}  $, then
$\Phi$ has a sequence of critical points $u_{k}\neq0$ such that $\Phi
(u_{k})\leq0$, $u_{k}\rightarrow0$.
\end{thm}

As pointed out in He--Wu \cite[Remark 2.5]{MR4162412}, in Theorem \ref{kt2},
instead of $\left(  PS\right)  $ condition, it suffices to assume $\left(
PS\right)  _{c}$ for $c\leq0$. That is, any sequence $\left\{  u_{n}\right\}
$ such that $\Phi^{\prime}(u_{n})\rightarrow0$ and $\Phi(u_{n})\rightarrow
c\leq0$, possesses a convergent subsequence.

We need the following lemma.

\begin{lem}
\label{kl1}If $u_{n}\rightharpoonup u$ in $E$, then%
\begin{equation}
\varliminf_{n\rightarrow\infty}\left[  \left(  \int\left\vert \nabla
u_{n}\right\vert ^{2}\right)  \int\nabla u_{n}\cdot\nabla(u_{n}-u)-\left(
\int\left\vert \nabla u\right\vert ^{2}\right)  \int\nabla u\cdot\nabla
(u_{n}-u)\right]  \geq0\text{,} \label{x}%
\end{equation}

\end{lem}

\begin{pf}
By direct computation we have%
\begin{align*}
& \left(  \int\left\vert \nabla u_{n}\right\vert ^{2}\right)
\int\nabla u_{n}\cdot\nabla(u_{n}-u)-\left(  \int\left\vert \nabla
u\right\vert ^{2}\right)  \int\nabla u\cdot\nabla(u_{n}-u)\\
&  =\left(  \int\left\vert \nabla u_{n}\right\vert ^{2}\right)  \int\left\vert
\nabla(u_{n}-u)\right\vert ^{2}+\left(  \int\left\vert \nabla u_{n}\right\vert
^{2}-\int\left\vert \nabla u\right\vert ^{2}\right)  \int\nabla u\cdot
\nabla(u_{n}-u)\\
&  \geq\left(  \int\left\vert \nabla u_{n}\right\vert ^{2}-\int\left\vert
\nabla u\right\vert ^{2}\right)  \int\nabla u\cdot\nabla(u_{n}-u)\text{.}%
\end{align*}
Since $u_{n}\rightharpoonup u$ in $E$, the right hand side goes to zero. The
desired result follows from taking lower limit on both sides of the above inequality.
\end{pf}

Now, we are ready to prove Theorem \ref{kt1}.

\begin{pf}
[\indent Proof of Theorem \ref{kt1}]Let $\phi:[0,\infty)\rightarrow\mathbb{R}$ be a decreasing
$C^{\infty}$-function such that $\left\vert \phi^{\prime}(t)\right\vert \leq
2$,%
\[
\phi(t)=1\text{\quad for }t\in\left[  0,1\right]  \text{,\qquad}%
\phi(t)=0\text{\quad for }t\geq2\text{.}%
\]
We consider the following truncated functional $I:E\rightarrow\mathbb{R}$,
which is a modification of the truncated functional used in \cite{MR4162412},
\begin{equation}
I(u)=\frac{1}{2}\left\Vert u\right\Vert ^{2}-\frac{1}{2}\left(  \left\Vert
u^{\ast}\right\Vert ^{2}+2\int F(u)\right)  \phi(\left\Vert u\right\Vert
^{2})+\frac{1}{4}\left(  \int\left\vert \nabla u\right\vert ^{2}\right)
^{2}\text{,} \label{I}%
\end{equation}
where $u^{\ast}=u^{-}+u^{0}\in E^{-}\oplus E^{0}$. The derivative $I^{\prime}$
is given by%
\begin{align}
\langle I^{\prime}(u),v\rangle  &  =\left[  1-\left(  \left\Vert
u^{\ast}\right\Vert ^{2}+2\int F(u)\right)  \phi^{\prime}(\left\Vert
u\right\Vert ^{2})\right]  \left(  u,v\right) \nonumber\\
&  \qquad-\left[  \left(  u^{\ast},v^{\ast}\right)  +\int f(u)v\right]
\phi(\left\Vert u\right\Vert ^{2})+\left(  \int\left\vert \nabla u\right\vert
^{2}\right)  \int\nabla u\cdot\nabla v \label{k5}%
\end{align}
for $u,v\in E$.

We will apply Theorem \ref{kt2} to $I$ to get a sequence of critical points
$\left\{  u_{k}\right\}  $ for $I$ such that%
\[
I(u_{k})\leq0\text{,\qquad}u_{k}\rightarrow0\text{.}%
\]
Since $I(u)=\Phi(u)$ for $\left\Vert u\right\Vert \leq1$, we see that for
large $k$ all the $u_{k}$ are critical points of $\Phi$ and Theorem \ref{kt1}
is proved.

Obviously $I$ is even. If $\left\Vert u\right\Vert \geq2$, then $\phi
(\left\Vert u\right\Vert ^{2})=0$. Hence%
\begin{align*}
I(u)  &  =\frac{1}{2}\left\Vert u\right\Vert ^{2}+\frac{1}{4}\left(
\int\left\vert \nabla u\right\vert ^{2}\right)  ^{2}\\
&  \geq\frac{1}{2}\left\Vert u\right\Vert ^{2}\rightarrow+\infty\text{,\qquad
as }\left\Vert u\right\Vert \rightarrow\infty\text{.}%
\end{align*}
This means that $I$ is coercive.

To verify $\left(  PS\right)  _{c}$ for $c\leq0$, let $\left\{  u_{n}\right\}
$ be a sequence in $E$ such that $
I(u_{n})\rightarrow c\leq0$, $I^{\prime}(u_{n})\rightarrow0$. 
Since $I$ is coercive, $\left\{  u_{n}\right\}  $ is bounded in $E$. Up to a
subsequence, we may assume that $u_{n}\rightharpoonup u$ in $E$. Then%
\[
-\left(  \left\Vert u_{n}^{\ast}\right\Vert ^{2}+2\int F(u_{n})\right)
\phi(\left\Vert u_{n}\right\Vert ^{2})=2I(u_{n})-\left\Vert u_{n}\right\Vert
^{2}-\frac{1}{2}\left(  \int\left\vert \nabla u_{n}\right\vert ^{2}\right)
^{2}\leq0\text{.}%
\]
Hence%
\begin{equation}
\left\Vert u_{n}^{\ast}\right\Vert ^{2}+2\int F(u_{n})\geq0\text{.} \label{kp}%
\end{equation}
Because $\phi^{\prime}(\left\Vert u_{n}\right\Vert ^{2})\leq0$ and%
\[
\varliminf_{n\rightarrow\infty}\left(  u_{n},u_{n}-u\right)  =\varliminf
_{n\rightarrow\infty}\left\Vert u_{n}\right\Vert ^{2}-\left\Vert u\right\Vert
^{2}\geq0\text{,}%
\]
up to a further subsequence we may assume
\begin{equation}
\left(  \left\Vert u_{n}^{\ast}\right\Vert ^{2}+2\int F(u_{n})\right)
\phi^{\prime}(\left\Vert u_{n}\right\Vert ^{2})\left(  u_{n},u_{n}-u\right)
\longrightarrow\alpha\leq0\text{,} \label{k7}%
\end{equation}
note here that by the boundedness of $\left\{  u_{n}\right\}  $, the
coefficient of $\left(  u_{n},u_{n}-u\right)  $ is bounded.

Thanks to Lemma \ref{kl1}, we may also assume%
\begin{equation}
\left(  \int\left\vert \nabla u_{n}\right\vert ^{2}\right)  \int\nabla
u_{n}\cdot\nabla(u_{n}-u)-\left(  \int\left\vert \nabla u\right\vert
^{2}\right)  \int\nabla u\cdot\nabla(u_{n}-u)\longrightarrow\beta\geq0\text{.}
\label{k6}%
\end{equation}
From the subcritical assumption $\left(  f_{1}\right)  $ and the compact
embedding $E\hookrightarrow L^{2}(\Omega)$, it is well known that%
\begin{equation}
\int f(u_{n})\left(  u_{n}-u\right)  \rightarrow0\text{,\qquad}\int
f(u)\left(  u_{n}-u\right)  \rightarrow0\text{.} \label{k8}%
\end{equation}
Finally, because $\dim(E^{-}\oplus E^{0})<\infty$, \ we also have%
\begin{equation}
\left(  u_{n}^{\ast},u_{n}^{\ast}-u^{\ast}\right)  \rightarrow0\text{,\qquad
}\left(  u^{\ast},u_{n}^{\ast}-u^{\ast}\right)  \rightarrow0\text{.}
\label{k9}%
\end{equation}
Computing $\left\langle I^{\prime}(u_{n}),u_{n}-u\right\rangle $ and
$\left\langle I^{\prime}(u),u_{n}-u\right\rangle $ via (\ref{k5}) then
subtracting the results, we deduce from (\ref{k7}), (\ref{k6}), (\ref{k8}) and
(\ref{k9}) that
\begin{align}
&\left\Vert u_{n}-u\right\Vert ^{2} =\langle I^{\prime
}(u_{n})-I^{\prime}(u),u_{n}-u\rangle \nonumber\\
&  \qquad\qquad+\left(  \left\Vert u_{n}^{\ast}\right\Vert ^{2}+2\int F(u_{n}%
)\right)  \phi^{\prime}(\left\Vert u_{n}\right\Vert ^{2})\left(  u_{n}%
,u_{n}-u\right) \nonumber\\
&  \qquad\qquad-\left(  \left\Vert u^{\ast}\right\Vert ^{2}+2\int F(u)\right)
\phi^{\prime}(\left\Vert u\right\Vert ^{2})\left(  u,u_{n}-u\right)
\nonumber\\
&  \qquad\qquad+\left[  \left(  u_{n}^{\ast},u_{n}^{\ast}-u^{\ast}\right)  +\int
f(u_{n})\left(  u_{n}-u\right)  \right]  \phi(\left\Vert u_{n}\right\Vert
^{2})\nonumber\\
&  \qquad\qquad-\left[  \left(  u^{\ast},u_{n}^{\ast}-u^{\ast}\right)  +\int
f(u)\left(  u_{n}-u\right)  \right]  \phi(\left\Vert u\right\Vert
^{2})\nonumber\\
&  \qquad\qquad-\left(  \int\left\vert \nabla u_{n}\right\vert ^{2}\right)
\int\nabla u_{n}\cdot\nabla(u_{n}-u)+\left(  \int\left\vert \nabla
u\right\vert ^{2}\right)  \int\nabla u\cdot\nabla(u_{n}-u)\nonumber\\
&  =\left[  o(1)+\alpha-\beta\right]  \rightarrow\left(  \alpha-\beta\right)
\leq0\text{.} \label{J}%
\end{align}
It follows that $u_{n}\rightarrow u$ in $E$ and $I$ satisfies $\left(
PS\right)  _{c}$ for $c\leq0$.

Finally, for $k\in\mathbb{N}$, let $X_{k}$ be an arbitrary $k$-dimensional
subspace of $E$. There is $\Lambda_{k}>0$ such that%
\[
\left\vert u\right\vert _{2}^{2}\geq\Lambda_{k}\left\Vert u\right\Vert
^{2}\text{\qquad for }u\in X_{k}\text{.}%
\]
There is also a constant $\eta>0$ such that for all $u\in E$ we have%
\[
\int\left\vert \nabla u\right\vert ^{2}\leq\eta\left\Vert u\right\Vert
^{2}\text{.}%
\]
From $\left(  f_{2}\right)  $, there is $\delta>0$ such that%
\begin{align}
F(t)\geq\frac{1+\eta^{2}}{\Lambda_{k}}t^{2}\text{\qquad for }t\in\left(
-\delta,\delta\right)  \text{.}\label{F}%
\end{align}
Take $\rho_{k}\in\left(  0,1\right)  $ such that if $u\in X_{k}$, $\left\Vert
u\right\Vert =\rho_{k}$, then $\left\vert u\right\vert _{\infty}<\delta$. For
$u\in X_{k}\cap S_{\rho_{k}}$ we have $|u(x)|\le\delta$ for all $x\in\Omega$.
Hence by (\ref{F}),
\begin{align*}
I(u)=\Phi(u)  &  =\frac{1}{2}\left(  \Vert u^{+}\Vert^{2}-\left\Vert
u^{-}\right\Vert ^{2}\right)  +\frac{1}{4}\left(  \int\left\vert \nabla
u\right\vert ^{2}\right)  ^{2}-\int F(u)\\
&  \leq\left\Vert u\right\Vert ^{2}+\frac{\eta^{2}}{4}\left\Vert u\right\Vert
^{4}-\frac{1+\eta^{2}}{\Lambda_{k}}\int u^{2}\\
&  \leq\frac{\eta^{2}}{4}\rho_{k}^{4}-\eta^{2}\rho_{k}^{2}\leq-\frac{3\eta
^{2}}{4}\rho_{k}^{2}\text{.}%
\end{align*}
Thus%
\[
\sup_{X_{k}\cap S_{\rho_{k}}}I\leq-\frac{3\eta^{2}}{4}\rho_{k}^{2}<0\text{.}%
\]
Now, by Theorem \ref{kt2}, $I$ has a sequence of critical points $\left\{
u_{k}\right\}  $ such that $u_{k}\rightarrow0$ in $E$. For some $k_{0}$, if
$k\geq k_{0}$ then $\left\Vert u_{k}\right\Vert <1$ and $u_{k}$ is a critical
point of $\Phi$. Hence $\Phi$ has a sequence of critical points $\{u_{k}%
\}_{k\ge k_{0}}$ converging to zero.
\end{pf}

\section{Proof of Theorem \ref{2k}}

Given $u\in E$, let $\phi_{u}$ be the solution of the second equation in the
system (\ref{kk}). It is well known that if $u\in E$ is a critical point of
$\Phi:E\rightarrow\mathbb{R}$,%
\begin{align*}
\Phi(u)  &  =\frac{1}{2}\int\left(  \left\vert \nabla u\right\vert
^{2}+V(x)u^{2}\right)  +\frac{1}{4}\int\phi_{u}u^{2}-\int F(u)\\
&  =\frac{1}{2}\left(  \Vert u^{+}\Vert^{2}-\left\Vert u^{-}\right\Vert
^{2}\right)  +\frac{1}{4}\int\phi_{u}u^{2}-\int F(u)\text{,}%
\end{align*}
then $\left(  u,\phi_{u}\right)  $ is a solution of (\ref{kk}), this idea was initiated from Benci \emph{et al.} \cite{MR1659454}. Similar to
(\ref{I}) we consider a truncated functional $I:E\rightarrow\mathbb{R}$
\[
I(u)=\frac{1}{2}\left\Vert u\right\Vert ^{2}-\frac{1}{2}\left(  \left\Vert
u^{\ast}\right\Vert ^{2}+2\int F(u)\right)  \phi(\left\Vert u\right\Vert
^{2})+\frac{1}{4}\int\phi_{u}u^{2}\text{.}%
\]
Then $I$ is an even coercive functional with $I(0)=0$. Similar to the last
section, using $\left(  f_{2}\right)  $, for $k$-dimensional subspace $X_{k}$
there is $\rho_{k}>0$ such that (\ref{k4}) holds.

To verify the $\left(  PS\right)  _{c}$ condition with $c\le0$ for $I$, we
need the following analogue of Lemma \ref{kl1}.

\begin{lem}
\label{kl2}If $u_{n}\rightharpoonup u$ in $E$, then%
\begin{equation}
\lim_{n\rightarrow\infty}\left(  \int\phi_{u_{n}}u_{n}\left(  u_{n}-u\right)
-\int\phi_{u}u\left(  u_{n}-u\right)  \right)  =0\text{.} \label{e}%
\end{equation}

\end{lem}

\begin{pf}
It is well known that $\phi_{u}$ is obtained from applying Riesz lemma to the
functional $\ell_{u}:v\mapsto\int u^{2}v$ on $E$. Thus%
\begin{align}
\left\Vert \phi_{u}\right\Vert  &  =\left\Vert \ell_{u}\right\Vert
=\sup_{\left\Vert v\right\Vert =1}\left\vert \int u^{2}v\right\vert
\nonumber\\
&  \leq\sup_{\left\Vert v\right\Vert =1}\left(\vert u^{2}\vert _{3}\left\vert v\right\vert _{3/2}\right)%
=\left\vert u\right\vert _{6}^{2}\sup_{\left\Vert v\right\Vert =1}\left\vert v\right\vert _{3/2}\leq
C\left\Vert u\right\Vert ^{2}\text{.} \label{f}%
\end{align}
Since $\left\{  u_{n}\right\}  $ is bounded, we know that $\left\{
\phi_{u_{n}}\right\}  $ is also bounded in $E$. By the compactness of the embedding $E\hookrightarrow L^{12/5}(\Omega)$, up to a subsequence we have
$u_{n}\rightarrow u$ in $L^{12/5}(\Omega)$. Hence%
\[
\left\vert \int\phi_{u_{n}}u_{n}\left(  u_{n}-u\right)  \right\vert
\leq\left\vert \phi_{u_{n}}\right\vert _{6}\left\vert u_{n}\right\vert
_{12/5}\left\vert u_{n}-u\right\vert _{12/5}\rightarrow0\text{,}
\]
because $\left\{  \phi_{u_{n}}\right\}  $ and $\left\{  u_{n}\right\}  $ are
bounded in $L^{6}(\Omega)$ and $L^{12/5}(\Omega)$, respectively. Similarly,
the second integral in (\ref{e}) vanishes as $n\rightarrow\infty$.
\end{pf}

Let $\left\{  u_{n}\right\}  $ be a $\left(  PS\right)  _{c}$ sequence of
$\Phi$ with $c\leq0$. It is easy to see that (\ref{kp}) still holds in current
situation, thus we have (\ref{k7}). Using (\ref{k7}), (\ref{k8}), (\ref{k9}),
and Lemma \ref{kl2} we have an analogue of (\ref{J})%
\[
\left\Vert u_{n}-u\right\Vert^2 \rightarrow\alpha\leq0\text{.}%
\]
Thus $u_{n}\rightarrow u$ in $E$ and $\left(  PS\right)  _{c}$ is verified.
Applying Theorem \ref{kt2}, $I$ has a sequence of critical points
$u_{k}\rightarrow0$. Since $I(u)=\Phi(u)$ for $\left\Vert u\right\Vert \leq1$,
for large $k$, $u_{k}$ is critical point of $\Phi$. Thus $\Phi$ has a sequence
of critical points $u_{k}\rightarrow0$ in $E$. From (\ref{f}) we have $\phi_{u_{k}%
}\rightarrow0$ in $E$. Thus (\ref{kk}) has a sequence of solutions $(u_{k}%
,\phi_{u_{k}})\rightarrow\left(  0,0\right)  $ in $E\times E$.

\end{document}